\documentclass[12pt,centertags,oneside]{article}
\usepackage{amsmath,amstext,amsthm,amscd,typearea}
\usepackage{amssymb}
\usepackage{a4wide}
\usepackage[mathscr]{eucal}
\usepackage{mathrsfs}
\usepackage{typearea}
\usepackage{charter}
\usepackage{pdfsync}
\usepackage{url}

\usepackage{xcolor}

\usepackage[a4paper,width=16.2cm,top=3cm,bottom=3cm]{geometry}

\numberwithin{equation}{section}

\newtheorem{theorem}{Theorem}[section]
\newtheorem{definition}[theorem]{Definition}
\newtheorem{proposition}[theorem]{Proposition}
\newtheorem{corollary}[theorem]{Corollary}
\newtheorem{lemma}[theorem]{Lemma}
\newtheorem{remark}[theorem]{Remark}

\newcommand{\supp}{{\rm Supp}}

\newcommand{\const}{\mathop{\mathrm{const}}\nolimits}

\newcommand{\ddc}{dd^c}
\newcommand{\dc}{d^c}

\newcommand{\codim}{{\rm codim\ \!}}

\newcommand{\capK}{\text{cap}}

\newcommand{\C}{\mathbb{C}}

\newcommand{\N}{\mathbb{N}}

\newcommand{\R}{\mathbb{R}}
\renewcommand\P{\mathbb{P}}

\title{\bf  Density currents and relative non-pluripolar products}
\providecommand{\keywords}[1]{\textbf{\textit{Keywords:}} #1}
\providecommand{\subject}[1]{\textbf{\textit{Mathematics Subject Classification 2010:}} #1}

\author{Duc-Viet Vu}
\newcommand{\Addresses}{{
		\bigskip
		\footnotesize
		\textsc{Duc-Viet Vu, University of Cologne, Division of Mathematics, Department of Mathematics and Computer Science, Weyertal 86-90, 50931, K\"oln,  Germany  \& Thang Long Institute of Mathematics and Applied Sciences, Hanoi, Vietnam}
		\noindent
		\par\nopagebreak
		\noindent
		\textit{E-mail address}: \texttt{vuduc@math.uni-koeln.de}	
}}

\date{\today}
\begin{document}
\maketitle
\begin{abstract} Let $X$ be a compact K\"ahler manifold of dimension $n$. For $1 \le j \le m$,  let $T_j$ be a closed positive $(1,1)$-current  on $X$ whose cohomology class is K\"ahler. Let $T$ be a closed positive current of bi-degree $(p,p)$ with $m+p \le n$. We prove that if $T_1,\ldots, T_m$ are of $T$-relative full mass intersection, then the $T$-relative non-pluripolar product of $T_1,\ldots, T_m$  is equal to the Dinh-Sibony product of $T_1,\ldots, T_m,T$.  
\end{abstract}
\noindent
, {full mass intersection}, {density current}, {Dinh-Sibony product}.
\\

\noindent
\subject{32Uxx}, {32Q15}.

%\tableofcontents

%%%%%%%%%%%%%%%%%%%%%%%%%%%%%
%%%%%%%%%%%%%%%%%%%%%%%%%%%%%%%%%%

\section{Introduction}
%Let $\omega$ be a K\"ahler form on $X$.

The purpose of this paper is to compare the notion of relative non-pluripolar products of currents given in \cite{BT_fine_87,BEGZ,GZ-weighted,Viet-generalized-nonpluri} and that of density currents introduced in \cite{Dinh_Sibony_density}. This can be seen as a follow-up of our previous works \cite{VietTuanLucas,Viet_Lucas}.

Let $X$ be a compact K\"ahler manifold of dimension $n$ and $T_1, \ldots, T_m$ closed positive $(1,1)$-currents on $X$. Let $T$ be a closed positive $(p,p)$-current on $X$ with $m+p \le n$. The  \emph{$T$-relative non-pluripolar product} $\langle \wedge_{j=1}^m T_j \dot{\wedge} T \rangle$ of $T_1, \ldots, T_m$, which was introduced in \cite{Viet-generalized-nonpluri}, is a closed positive current of bi-degree $(m+p,m+p)$. We can view it as a sort of intersection of $T_1,\ldots,T_m,T$.  When $T \equiv 1$, we recover  the usual non-pluripolar product of $T_1, \ldots,T_m$ given in \cite{BT_fine_87,BEGZ,GZ-weighted}.

To continue, we need to recall some terminologies and notations.  For closed positive currents $S$, we denote by $\{S\}$ the cohomology class of $S$.   For cohomology classes $\alpha,\beta$ in $H^{p,p}(X,\R)$, we write $\alpha \le \beta$ if $(\beta-\alpha)$ can be represented by a closed positive current. Given two closed positive $(1,1)$-currents $S,S'$,   we say that $S'$ is \emph{less singular than} $S$ if every (global) potential of $S$ is bounded from above by a potential of $S'$. A currents $S'$ in a pseudoeffective  $(1,1)$-class $\alpha$  is said to have \emph{minimal singularities} if for every current $S$ in $\alpha$, then $S'$ is less singular than $S$. Such a current always exists by taking a suitable envelop of quasi-plurisubharmonic functions. We refer to \cite{Demailly_analyticmethod} for details.

A crucial feature of relative non-pluripolar products is that they satisfy the following monotonicity: for every closed positive $(1,1)$-current $T'_j$ in the same cohomology class of $T_j$ such that $T_j'$ is less singular than $T_j$ for $1 \le j \le m$, we have  
\begin{align}\label{ine-mono-relativenopluri}
\{\langle T_1 \wedge \cdots \wedge T_m \dot{\wedge} T\rangle \} \le \{\langle T_1' \wedge \cdots \wedge T_m' \dot{\wedge} T \rangle \}.
\end{align}
We refer to \cite[Theorem 1.1]{Viet-generalized-nonpluri} for a proof; see also \cite{BEGZ,Lu-Darvas-DiNezza-mono,WittNystrom-mono} for weaker versions. % and  \cite[Remark 2.2]{Viet-lelong-full-mass} for a variant of (\ref{ine-mono-relativenopluri}).  
Recall that $T_1,\ldots, T_m$ are said to be  of \emph{$T$-relative full mass intersection} if the equality  in (\ref{ine-mono-relativenopluri}) occurs for $T'_j$ to be a current with minimal singularities in the cohomology class of $T_j$. This is independent of the choice of currents with minimal singularities.

We have another general way to look at the intersection of $T_1, \ldots, T_m,T$ provided by the theory of density currents coined in \cite{Dinh_Sibony_density}. The last theory has found deep applications in complex dynamics and holomorphic foliations; for example, see \cite{DS_unique_ergodicity,DNS_foliation,VA-lyapunov,DNV,DNT_equi,DVT_growth_periodic,lucas_foliation}. We refer to the next section for basic facts about density currents.  In nice situations, the theory of density currents gives rise to a new notion of intersection of currents which is called \emph{the Dinh-Sibony product}. It was proved in \cite{DNV,VietTuanLucas,Viet_Lucas} that the last product covers classical intersections of currents. The following main result of this paper shows that this is also the case for currents of full mass intersection in K\"ahler classes.    

\begin{theorem}\label{the-main1} Assume that the cohomology class of $T_j$ is K\"ahler for $1 \le j \le m$. Then,   $T_1, \ldots, T_m$ are of $T$-relative full mass intersection  if and only if    we have 
$$\langle T_1 \wedge \cdots \wedge T_m \dot{\wedge} T\rangle=T_1 \curlywedge  \cdots \curlywedge  T_m \curlywedge T.$$
\end{theorem}

Note that the above result is no longer true if the classes of $T_j$'s are not K\"ahler.  This can be seen by using the fact that currents with minimal singularities in a non-nef  pseudoeffective $(1,1)$-class must have strictly positive Lelong number somewhere; see \cite[Corollary 6.4]{Demailly_regula_11current}. However, it is still plausible to expect that Theorem \ref{the-main1} still holds if $T\equiv 1$ and the classes of $T_j$'s are big nef. Unfortunately, we don't know how to prove it.

The uniform quasi-continuity of bounded psh functions proved in \cite{Viet-generalized-nonpluri} will play a key role in our proof of Theorem \ref{the-main1}; see Theorem \ref{th-capconvergedecreasing} below. This property allows us to prove that  in general, the relative non-pluripolar product  $\langle \wedge_{j=1}^m T_j \dot{\wedge} T \rangle$ is always bounded from above by density currents associated to $T_1,\ldots, T_m,T$.  We note that the case where $T\equiv 1$ was treated in \cite{Viet_Lucas}.  \\

\noindent
\textbf{Acknowledgments.} This research  is supported by a postdoctoral fellowship of the Alexander von Humboldt Foundation. \\

\section{Preliminaries on density currents}

In this section, we recall some basic properties of density currents. The last notion was introduced in \cite{Dinh_Sibony_density}. Most of materials below is taken from the last paper. Readers can also consult \cite{Vu_density-nonkahler} for some simplifications. 

Let $X$ be a complex K\"ahler manifold of dimension $n$ and $V$ a smooth complex submanifold of $X$ of dimension $l.$  Let $T$ be a closed positive $(p,p)$-current on $X,$ where $0 \le p \le n.$  Denote by $\pi: E\to V$ the normal bundle of $V$ in $X$ and $\overline E:= \P(E \oplus \C)$ the projective compactification of $E.$ By abuse of notation, we also use $\pi$ to denote the natural projection from $\overline E$ to $V$. 

Let $U$ be an open subset of $X$ with $U \cap V \not = \varnothing.$  Let $\tau$ be  a smooth diffeomorphism  from $U$ to an open neighborhood of $V\cap U$ in $E$ which is identity on $V\cap U$ such that  the restriction of its differential $d\tau$ to $E|_{V \cap U}$ is identity.  Such a map is called \emph{an admissible map}.  Note that in \cite{Dinh_Sibony_density}, to define an admissible map,  it is required furthermore that $d\tau$ is $\C$-linear at every point of $V$. This difference doesn't affect what follows.  When $U$ is a small enough tubular neighborhood of $V,$ there always exists an admissible map $\tau$ by \cite[Lemma 4.2]{Dinh_Sibony_density}. In general, $\tau$ is not holomorphic.  When $U$ is a small enough local chart, we can choose a holomorphic admissible map by using suitable holomorphic coordinates on $U$.   For $\lambda \in \C^*,$ let $A_\lambda: E \to E$ be the multiplication by $\lambda$ on fibers of $E.$

Let $\tau$ be an admissible map defined on a tubular neighborhood of $V$. It was proved in \cite[Theorem 4.6]{Dinh_Sibony_density} that the family $(A_\lambda)_* \tau_* T$ is of mass uniformly bounded in $\lambda$ and every limit current is a closed positive current on $\overline E$. We call these limits \emph{tangent currents to $T$ along $V$}.  Moreover, if 
$$S=\lim_{k\to \infty} (A_{\lambda_k})_* \tau_* T$$ for some sequence $(\lambda_k)_k$ converging to $\infty$, then  for every open subset $U$ of $X$ and  every admissible map $\tau': U' \to E$ , we also have  
$$S=\lim_{k\to \infty} (A_{\lambda_k})_* \tau'_* T.$$
This is equivalent to saying that tangent currents are independent of the choice of the admissible map $\tau$.  By this reason, the sequence $(\lambda_k)_k$ is called \emph{a defining sequence} of $T_\infty.$  
In practice, we usually choose $\tau_j$ to be a change of coordinates. Another crucial property of tangent currents to $T$ along $V$ is that although in general they are no unique, their cohomology classes are unique.   We call the last class  \emph{the total tangent class of $T$ along $V$}.

\begin{definition} (\cite{Dinh_Sibony_density}) Let $F$ be a complex manifold and $\pi_F: F \to V$ a holomorphic submersion. Let $S$  be  a positive current $S$ of bi-degree $(p,p)$  on $F$. \emph{The h-dimension} of $S$ with respect to $\pi_F$ is the biggest integer $q$ such that $S \wedge \pi_F^* \theta^q \not =0$ for some Hermitian metric $\theta$ on $V$.  
\end{definition}

%The h-dimension of $S$ is in $[\max\{l- p,0\}, \min\{\dim F -p,l\}]$. Recall that by  \cite[Lemma 3.2]{Dinh_Sibony_density}, given two nonnegative integers $m,q$, the h-dimension of $R$ is strictly less than $\max\{m,q\}$ if and only if $R \wedge \pi^* \theta=0$ for every smooth $(m,q)$-form $\theta$.  
We have the following description of currents with minimal h-dimension. 

\begin{lemma} \label{le-minimalhdimension} (\cite[Lemma 3.4]{Dinh_Sibony_density})  Let $\pi_F: F \to V$ be a submersion. Let $S$ be a closed positive  current of bi-degree $(p,p)$ on $F$ of h-dimension $(l -p)$ with respect to $\pi_F$. Then $S= \pi^* S'$ for some closed positive current $S'$ on $V$. 
\end{lemma}

Let $m\in \N^*$. Let $T_j$ be a closed positive current  of bi-degree $(p_j, p_j)$ for $1 \le j \le m$ on $X$ and  $T_1 \otimes \cdots \otimes T_m$ the tensor current of $T_1, \ldots, T_m$ which is a current on $X^m.$ %The current $T_1 \otimes \cdots \otimes T_m$ has mass on the diagonal $\Delta_m$ of $X^m$ if and only if $p_j=n$ (i.e, $T_j$'s are measures) and there exists $a \in X$ such that $T_j$ has Dirac mass at $a$ for every $1 \le j \le m$.  We will \emph{never} encounter the last situation in this paper.  
A \emph{density current} associated to $T_1, \ldots,  T_m$ is a tangent current to $\otimes_{j=1}^m T_j$ along the diagonal $\Delta_m$ of $X^m.$ Let $\pi_m: E_m \to \Delta$ be the normal bundle of $\Delta_m$ in $X^m$. The unique cohomology class of density currents associated to $T_1,\ldots,T_m$ is called \emph{the total density class of $T_1, \ldots, T_m$}.  When $m=2$, $T_2 =[V]$, the density currents of $T_1, T_2$ are naturally identified with the   tangent currents to $T_1$ along $V$ (see \cite[Lemma 2.3]{Vu_density-nonkahler}).  

  %The h-dimension of the total density class of   $T_1, \ldots T_m$ is called \emph{the density h-dimension} of  $T_1, \ldots T_m$.  

\begin{definition}  \label{def-DSproduct} We say that the \emph{Dinh-Sibony product} $T_1 \curlywedge \cdots \curlywedge T_m$ of $T_1, \ldots, T_m$ is well-defined  if $\sum_{j=1}^m p_j \le n$ and  there is only one density current associated to $T_1, \ldots, T_m$ and this current is  the pull-back by $\pi_m$ of a current $S$ on $\Delta_m$. We define $T_1 \curlywedge \cdots \curlywedge T_m$ to be $S$. 
\end{definition}

The notion of Dinh-Sibony products generalizes well-known notions of intersections of currents. We refer to \cite{DNV,Dinh_Sibony_density,VietTuanLucas} for details. We only cite here the following particular case of main results in \cite{VietTuanLucas} (see Remark \ref{re-noncompact} below for a definition of density currents in a general complex manifolds).
  
\begin{theorem}\label{th-density11agreewithclassical} Let $X$ be a complex manifold. Let $T_1, \ldots, T_m$ be closed positive currents of bi-degree $(1,1)$ on $X$ with locally bounded potentials and $T$ a closed positive current on $X$. Then the Dinh-Sibony product of $T_1, \ldots, T_m,T$ is well-defined and equal to the usual intersection $T_1 \wedge \cdots \wedge T_m \wedge T$. 
\end{theorem}

The following  results will be important later.
 
\begin{corollary} \label{cor-classDSproduct} (\cite[Page 546]{Dinh_Sibony_density}) Let $T_j$ be a closed positive current of bi-degree $(p_j,p_j)$ on $X$ for $1 \le j \le m$ such that  $\sum_{j=1}^m p_j \le n$.  Assume that the h-dimension of the total density class associated to $T_1, \ldots, T_m$ is minimal, \emph{i.e}, equal to $n- \sum_{j=1}^m p_j$. Then the total density class of  $T_1, \ldots, T_m$ is equal to $\pi_{m}^*(\wedge_{j=1}^m\{T_j\})$. 
\end{corollary}

We consider now a technical lemma which is useful later. Let $X$ be a complex manifold of dimension $n$.  Let $\omega$ be a Hermitian metric on $X$. Let $V$ be a smooth complex submanifold on $X$ and $\pi: E\to V$ the normal bundle of $V$ in $X$. Let $\sigma: \widehat X \to X$  be the blowup of $X$ along $V$ and $\widehat V$ the exceptional hypersurface. Let $W_1$ be an open subset of $X$ so that $W_1 \cap V \not  = \varnothing$  and $W_1 \cap V \Subset X$. 

Let $\widehat \omega_h$ be a Chern form of the dual of the  line bundle generated by $\widehat V$ on $\widehat X$ such that the restriction of $\widehat \omega_h$ to every fiber of the projection $\sigma|_{\widehat V}: \widehat V \to V$ is strictly positive. We have   
$$\widehat \omega: = c_0 \sigma^* \omega+ \widehat \omega_h >0$$
 on $\widehat W_1:= \sigma^{-1}(W_1)$ for some constant $c_0>0$.  Note that if $\omega$ is K\"ahler, then so is $\widehat \omega$ because $\widehat \omega_h$ is closed. Let $W_2$ be an open neighborhood of $V$ in $X$. Let $\varphi_V$ be a potential of $\sigma_* \widehat \omega_h$  with $\supp \varphi \Subset W_2$ if $\codim V \ge 2$ (this is always possible because $\sigma_* \widehat \omega_h$ is smooth outside $V$), otherwise let $\varphi_V$ be a smooth potential of $[V]$ compactly supported in $W_2$. 

Observe that there exists a constant $c>0$ such that $(\ddc \sigma^* \varphi_V - c [\widehat V])$ is a smooth form (see \cite[Page 509]{Dinh_Sibony_density}).  By rescaling $\widehat \omega_h$, we can assume that the constant $c$ in the above lemma is equal to $1$ or in other words, the current $(\ddc \sigma^* \varphi_V - [\widehat V])$ is smooth.   The function $\varphi_V$ plays a crucial role in estimating the mass of  tangent currents as showed in  \cite{Dinh_Sibony_density}.  

 Let $x:=(x',x'')$ be a local coordinate system on a local chart $U$ on $X$ with $V= \{x''=0\}$. Put $\gamma:= \varphi_V - \log \|x''\|$. Then $\sigma^* \gamma$ is smooth on  $\widehat U:= \sigma^{-1}(U)$. Particularly, $|\gamma|$ is uniformly bounded on $K \backslash V$ where $K$ is a compact of $U$. Thus  if $\|x''\| \le |\lambda|^{-1}$, then $e^{\varphi_V} \le e^A |\lambda|^{-1}$ on $K$, where $A:= \|\gamma\|_{L^\infty(K)}$.

\begin{lemma} \label{le-uocluongAlambda} (\cite[Lemma 2.12]{Dinh_Sibony_density}) Let $A>0$ be a constant.  There exists a sequence of smooth quasi-psh $(\varphi_{V, \lambda})_\lambda$ on $X$  depending only on $|\lambda|$ and decreasing to $\varphi_{V}$ as $|\lambda| \to \infty$ such that  $\supp \varphi_{V,\lambda} \subset \supp \varphi_V \Subset W_2$ and  for any local coordinate system $\big(U, x:=(x',x'')\big)$ on $W_1$ with $V= \{x''=0\}$ and $K \Subset U$, we have 
\begin{align}\label{ine-uonluongalambdax''V}
|\lambda|^2 \ddc \|x''\|^2 \le C_1 \big(\ddc \varphi_{V, \lambda} + C_1 \sigma_* \widehat \omega \big) 
\end{align}
on $\{e^{\varphi_V} \le e^A |\lambda|^{-1}\}\cap K$  for some constant $C_1>1$ independent of $\lambda$.
\end{lemma}

We conclude this section with the following remark.

\begin{remark} \label{re-noncompact} It is sometimes convenient to work with non-closed currents or work in the local setting. We provide here some terminology for later use. Assume $X$ is now not necessarily compact nor K\"ahler and $T$ a positive current on $X$. We say that $T$  \emph{ admits tangent currents along $V$} if for every local holomorphic admissible map $\tau: U \to E$, the family $ (A_{\lambda})_* \tau_* T$ is of uniformly bounded mass on compact subsets of $E|_{U \cap V}$. A \emph{tangent current}  $T_\infty$ to $T$ along $V$ is a current on $E$ such that there are a sequence $(\lambda_k)_k \subset \C^*$ converging to $\infty$ and a collection of holomorphic admissible maps $\tau_j: U_j \to E$ for $j \in J$ satisfying that  $V \subset \bigcup_{j \in J} U_j$ and 
$$T_\infty:= \lim_{k \to \infty} (A_{\lambda_k})_* (\tau_j)_* T$$
on $\pi^{-1}(U_j \cap V)$ for every $j \in J.$ In this setting, it is still true that tangent currents are independent of the choices of admissible maps; see \cite[Proposition 2.5]{Viet_Lucas} (the proof there requires only the positivity of  currents). 
%   have to consider limit currents of the family $(A_\lambda)_* \tau_*(\psi T)$, where $\psi$ is a Borel bounded function on $X$. Any limit current of this family is also called \emph{a tangent current to $ \psi T$ along $V$}. We will see an example of this kind of currents in Section \ref{subsec-phiV}. 
\end{remark}

\section{Proof of the main result}

To study the relation between relative non-pluripolar products and density currents, we need to digress for a moment to  recall some facts about quasi-continuity of bounded psh functions.

Let $U$ be an open subset of $\C^n$.   %We recall some notations and results from \cite{Kolodziej05}.  If 
Let $K$ be a Borel subset of  $U$.  %By  \cite{Bedford_Taylor_82}, the  \emph{capacity} $\capK(K,U)$ of $K$ in $U$ is given by $$\capK(K,U):= \sup \bigg\{ \int_K (\ddc u)^n:  u \text{ is psh  on } U \text{ and } 0 \le u \le 1\bigg\}.$$
For a  closed positive current $T$ of bi-dimension $(m,m)$ on $U$ ($0 \le m \le n$), recall that
$$\capK_T(K,U):= \sup \bigg\{ \int_K (\ddc u)^m \wedge T:  u \text{ is psh  on } U \text{ and } 0 \le u \le 1\bigg\},$$
see \cite{Bedford_Taylor_82,Kolodziej05,Xing-continuity}.  Note that the above notion makes sense if $T$ is a smooth (not necessarily closed) positive form. 

\begin{lemma} \label{le-capTK} (\cite[Lemma 2.1]{Viet-generalized-nonpluri})  Let $A$ be a locally complete pluripolar set in $U$. Let $T$ be a closed positive current of bi-dimension $(m,m)$ on $U$. Assume that  $T$ has no mass on $A$. Then,  we have $\capK_T(A,U)=0$. 
\end{lemma}

Here is the uniform strong quasi-continuity of bounded p.s.h. functions which we will need to use. 

\begin{theorem}\label{th-capconvergedecreasing} (a particular case of \cite[Theorem 2.4]{Viet-generalized-nonpluri}) Let  $T_k = \ddc w_k\wedge S$, $T= \ddc w\wedge S$, where  $S$ is a closed positive current,  $w$ is a psh function locally integrable with respect to $S$ and $w_k$ is a psh function converging to $w$ in $L^1_{loc}$ as $k \to \infty$ so that $w_k \ge w$ for every $k$.    Then,  for every constant $\epsilon>0$,  there exists an open subset $U'$ of $U$ such that  $\capK_{T_k}(U',U)<\epsilon$ for every $k$ and the restriction of  $u$ to $U \backslash U'$ is continuous. 
\end{theorem}

As a direct consequence, we obtain the following.  

\begin{corollary} \label{cor-convergplurifine} Let $T_k$ be as in Theorem \ref{th-capconvergedecreasing}. Let $T'_k$ be either a closed positive current of the same bi-degree as $T_k$ or a (not necessarily closed) smooth positive form of the same bi-degree as $T_k$  such that $T'_k \le T_k$ for every $k$.   Let $R_k:= \ddc v_{1k} \wedge \cdots \wedge \ddc v_{m k} \wedge T_k$ and $R:=  \ddc v_{1} \wedge \cdots \wedge  \ddc v_{m} \wedge T$, where $v_{j k}, v_j$ are bounded  psh functions on $U$ such that $v_{jk}$ is bounded uniformly  in $k$.   Let $u$ be a bounded psh function on $U$ and $\chi$ a continuous function on $\R$.  Assume that $R_k \to R$ as $ k\to \infty$ on $U$.
 Then we have 
$$\chi(u) R_k \to \chi(u) R$$
 as $k \to \infty$. In particular, the last convergence holds when $T_k= T$ for every $k$ or $T_k = \ddc w_k\wedge S$, $T= \ddc w\wedge S$, where  $S$ is a closed positive current,  $w$ is a psh function locally integrable with respect to $S$ and $w_k$ is a psh function converging to $w$ in $L^1_{loc}$ as $k \to \infty$ so that $w_k \ge w$ for every $k$.   
\end{corollary}

\proof Since $T'_k \le T_k$, the uniform  strong quasi-continuity of bounded psh functions also holds for $(T'_k)_k$. The proof now goes verbatim  as that of  \cite[Corollary 2.5]{Viet-generalized-nonpluri}. 
\endproof

Let $X$ be a complex manifold of dimension $n$ and  $m\in \N$.   Let $\Delta_{m+1}$ be the diagonal of $X^{m+1}$ and $\pi_{m+1}: E_{m+1}\to \Delta_{m+1}$ the normal bundle of $\Delta_{m+1}$ in $X^{m+1}$.  Let $T$ be a closed positive current of bi-degree $(p,p)$ on $X$. Let  $T_j$ be closed positive current of bi-degree $(1,1)$ for $1 \le j \le m$.  Let $p_j: X^{m+1} \to X$ be the natural projection from $X^{m+1}$ to its $j^{th}$-component for $1 \le j \le m+1$.  %  In views of applications to generalized non-pluripolar products,  we \emph{always} assume that $p+ m \le n$ because otherwise the generalized non-pluripolar product of $T_1, \ldots, T_m,T$ is zero by a bi-degree reason. 

 In what follows, we will use terminologies from Remark \ref{re-noncompact}.  Here is  a key auxiliary result for us.
   
\begin{lemma} \label{le-hoitutoidensityboundedpsh} Let $\psi$ be a locally bounded quasi-psh function on $X^{m+1}$. Let $R:= \otimes_{j=1}^m T_j \otimes T$.   Assume that $T_1, \ldots, T_m$ have locally bounded potentials. Then, the unique tangent current to $\psi R$ along $\Delta_{m+1}$ is $\pi_{m+1}^* \big((\psi|_{\Delta_{m+1}})  \wedge_{j=1}^m T_j \wedge T\big)$. 
\end{lemma}

Recall that $\psi|_{\Delta_{m+1}}$ is the restriction of $\psi$ to $\Delta_{m+1}$. 

\proof By Theorem \ref{th-density11agreewithclassical}, the tangent current to $R$ along $\Delta_{m+1}$ is equal to $\pi_{m+1}^* (\wedge_{j=1}^m T_{j} \wedge T)$. Hence, the desired assertion is  clear if $\psi$ is continuous. The whole arguments below are to justify that we still have the same conclusion for every locally bounded psh function $\psi$. This will be done by using the uniform strong quasi-continuity of bounded psh functions mentioned above. 

Since the problem is local,   we can assume $X\Subset \C^n$, $T_j= \ddc u_j$ for some  psh function $-1 \le u_j \le 0$ on $X$ for $1 \le j \le m$ and $\psi$ is bounded on $X$. Let $Y:= X^{m+1}$ and $\sigma: \widehat Y \to Y$ the blow-up of $Y$ along $V:= \Delta_{m+1}$. Let $\widehat V$ be the exceptional hypersurface. Let $\pi: E \to V$ be the normal bundle of $V$ in $Y$ (we change a bit our notations here for convenience). %Let $\overline E:= \P(E \oplus \C)$,  $H_\infty:= \overline E \backslash E$ which is naturally isomorphic to $\P(E)$ and $\pi_{\P(E)}: H_\infty \to V$ the natural projection.  
 Let $\widehat E$ be the blowup of $E$ along $V$.  Note that $\widehat E$ is naturally identified with the normal bundle of $\widehat V$ in $\widehat Y$.  Let $\pi_{\widehat E}: \widehat E \to \widehat V$ be the natural projection. %We use the identity maps as an admissible map when working with density currents.  

Let $\widehat R$ be the strict transform of $R$ by $\sigma$.  Put $$Q:= \wedge_{j=1}^m T_j \wedge T.$$
We recall that there is an one-to-one correspondence between tangent currents to $\widehat R$ along $\widehat V$ and those to $R$ along $V$, see \cite[Lemma 4.7]{Dinh_Sibony_density}. Thus, the unique tangent current to $\widehat R$ along $\widehat V$ is  $\pi_{\widehat E}^* (\sigma|_{\widehat V})^* Q,$
here we implicitly identify $V$ with $X$.   As a consequence, for every admissible map $\widehat \tau$ from an open subset of $\widehat Y$ to $\widehat E$ and every continuous function $g$ on $\widehat Y$, we have
\begin{align}\label{limit-tangentcurrentR}
(A_\lambda)_* \widehat \tau_* (g \widehat R)=  \pi_{\widehat E}^*  \big( (g|_{\widehat V})  (\sigma|_{\widehat V})^* Q\big),
\end{align}
where $A_\lambda$ is the fiberwise multiplication by $\lambda$ in $\widehat E$. 
Observe that $p_j \circ \sigma$ is submersion for every $j$. One can see it by using standard local coordinates on $\widehat Y$.  Put 
$$\widehat u_{j}:= (p_j\circ \sigma)^* u_{j}$$
 for every $j,k$, and 
$$\widehat T:= (p_{m+1}\circ \sigma)^* T.$$  Note $\widehat u_j$ is bounded by hypothesis.   Let $\varphi_{\widehat V}, \varphi_{\widehat V, \lambda}$ be functions in Lemma \ref{le-uocluongAlambda} for $\widehat V,\widehat Y$ in place of $V,X$ respectively. Recall that $\varphi_{\widehat V, \lambda}$ decreases to $\varphi_{\widehat V}$ as $\lambda \to \infty$ and $\varphi_{\widehat V}$ is a potential of $\widehat V$. 

Since the fibers of  $p_{m+1} \circ \sigma$ intersect $\widehat V$ properly, we see that $\varphi_{\widehat V}$ is locally integrable with respect to  $\widehat T$. In particular the last current  has no mass on $\widehat V$. It follows that the current 
$$\wedge_{j=1}^m \ddc \widehat u_{j} \wedge \widehat T$$
 has no mass on $\widehat V$ by  Lemma \ref{le-capTK}. % \cite[Lemma 2.1]{Viet-generalized-nonpluri}.
On the other hand, since $\sigma$ is isomorphic outside $\widehat V$, we get 
$$\wedge_{j=1}^m \ddc \widehat u_{j} \wedge \widehat T= \widehat R$$
on $\widehat Y \backslash \widehat V$.   This combined with the fact that the both sides of the last equality have no mass on $\widehat V$ yields that 
$$\wedge_{j=1}^m \ddc \widehat u_{j} \wedge \widehat T= \widehat R$$
 on $\widehat Y$. Using the last equality,  the local integrability of $\varphi_{\widehat V}$ with respect to $\widehat T$,  we deduce that $\varphi_{\widehat V}$ is locally integrable with respect to  $\widehat R$. 
This together with  Theorem \ref{th-capconvergedecreasing} implies that we have a uniform strong quasi-continuity property for $\widehat \psi:=\sigma^* \psi$  with respect to the sequence of currents $(\ddc \varphi_{\widehat V, \lambda}+ \widehat \omega) \wedge \widehat T$ as $\lambda \to \infty$.  

Let $\widehat U$ be a local chart of $\widehat Y$ small enough such that $\widehat E$ is trivializable on $\widehat V \cap \widehat U$. Thus, we can consider $\widehat U$ as an open subset of $\widehat E$ (by shrinking $\widehat U$ a bit if necessary) and use the identity map as a holomorphic admissible map from $\widehat U$ to $\widehat E$.

Let $q:= \dim Y- (m+p)$ and let $\Phi$ be a smooth positive form of bi-degree $(q,q)$ with compact support in $\widehat U$. By Lemma \ref{le-uocluongAlambda}, we obtain
$$ A_\lambda^* \Phi \lesssim  (\ddc \varphi_{\widehat V, \lambda}+ \widehat \omega) \wedge \widehat \omega^{q-1}.$$
Using this and  Corollary \ref{cor-convergplurifine} and the above uniform quasi-continuity of $\widehat \psi$,  we get a uniform strong quasi-continuity of  $\widehat \psi$  with respect to the sequence of positive currents  $\widehat T \wedge A_\lambda^* \Phi$ as $\lambda \to \infty$.  In other words, for  every constant $\epsilon>0$, there exists an open subset $\widehat U'$ in $\widehat U$ such that 
\begin{align}\label{ine-capKTwedgealmbdaphi}
\capK_{\widehat T \wedge A_\lambda^* \Phi}(\widehat U', \widehat U) \le \epsilon
\end{align}
 for every $\lambda$ and $\widehat \psi$ is continuous on $\widehat U \backslash \widehat U'$. Furthermore, we can also choose an open subset $W$ of $X\approx V$ so that  $\psi|_V$ is continuous on $X \backslash W$ and 
\begin{align}\label{ine-capKTwedgealmbdaphi2}
\capK_{T}(W) \le \epsilon.
\end{align}
 We deduce that $\widehat \psi$ is continuous on 
$$B:= (\widehat U \backslash \widehat U')\cup \sigma^{-1}(V \backslash W)$$
 (we identified $V$ with $X$). Let $\psi'$ be a bounded continuous function on $\widehat U$ such that $\psi'= \widehat \psi$ on $B\cap \widehat U$.  The inequality (\ref{ine-capKTwedgealmbdaphi}) ensures that 
$$\|(A_\lambda)_* (\widehat \psi \widehat R)- (A_\lambda)_* (\psi' \widehat R) \| \lesssim \epsilon,$$
 where the mass is measured on compact subsets of $\widehat U$ and  (\ref{ine-capKTwedgealmbdaphi2}) yields that 
$$\| (\widehat \psi|_{\widehat V}- \psi'|_{\widehat V})  (\sigma|_{\widehat V})^* (\wedge_{j=1}^m T_{j} \wedge T)\| \lesssim \epsilon,$$
note here that 
$$\widehat \psi|_{\widehat V}=(\psi\circ \sigma)|_{\widehat V}=  (\sigma|_{\widehat V})^* (\psi|_V).$$
This combined with (\ref{limit-tangentcurrentR}) applied to $g= \psi'$  gives the desired assertion. The proof is finished.  
\endproof

%The techniques used in the last proof  will be applied systematically in the next parts of this paper.    

 Here is a crucial connection between density currents and relative non-pluripolar products. 

\begin{theorem} \label{the-phanbucuarestricteddenstyvanonpluri}  Assume that $T_1, \ldots, T_m, T$ admit density currents. Let $R_\infty$ be a density current associated to $T_1, \ldots, T_m,T$. Then we have 
\begin{align} \label{ine_TjS2}
\pi^*_{m+1} \langle \wedge_{j=1}^m T_j  \dot{\wedge} T \rangle \le R_\infty
\end{align} 
(recall that we identify $\Delta_{m+1}$ with $X$). 
\end{theorem}

Note that (\ref{ine_TjS2}) was proved in \cite[Theorem 3.1]{Viet_Lucas} provided that $T$ is a constant function. The proof below will clarify some arguments there.

\proof  Since the problem is local, we can assume $X\Subset \C^n$ and every current in question is defined in an open neighborhood of $\overline X$ and  $T_j= \ddc u_j$ for some negative psh function $u_j$ on $X$ for $1 \le j \le m$.  Write $T_j= \ddc u_j+ \theta_j$, where $\theta_j$ is a smooth form and $u_j$ is a $\theta_j$-psh on $X$.  
Let 
$$\psi:=  \sum_{j=1}^m p_j^* u_j, \quad \psi_k:= k^{-1}\max\{\psi, -k\}$$ for $k \in \N$. Put $$u_{jk}:= \max\{u_j,-k\}, \quad T_{jk}:=\ddc u_{jk}$$
for $k\in \N$. Define  
$$R:= \otimes_{j=1}^m T_{j} \otimes T, \quad R_k:= \otimes_{j=1}^m T_{jk} \otimes T.$$
By definition, $R_\infty$ is a tangent current of $R$ along $\Delta_{m+1}$. To simplify the notation (without loss of generality), we  assume that $R_\infty$ is the unique density current associated to $T_1,\ldots, T_m,T$.

Let $\omega$ be a K\"ahler form on $(\C^n)^{m+1}$ and $q:= m(n-1)+ n-p$. Observe that 
$$\psi_k= k^{-1} \max\{ \sum_{j=1}^m u_{jk}, -k\}.$$
 Let $\epsilon>0$ be a constant. Using the uniform strong quasi-continuity for $u_{jk}$ (Theorem \ref{th-capconvergedecreasing}), we see that there is  a continuous function $v_{jk}$ on $X$ such that 
$$\capK_{T_{j}}(v_{jk} \not = u_{jk}) \le \epsilon, \quad \capK_{T_{jr}}(v_{jk} \not = u_{jk}) \le \epsilon$$
 for every $r\in \N$. Put 
$$\tilde{\psi}_{k}:=  k^{-1} \max\{ \sum_{j=1}^m v_{jk}, -k\}$$
 which is a continuous function. Notice that 
$$\{ \tilde{\psi}_k \not = \psi_k\} \subset \cup_{j=1}^m p_j^* \{ v_{jk} \not = u_{jk} \}.$$
 Thus 
$$\int_{X^{m+1}} |\psi_k - \tilde{\psi}_k| R_r \wedge \omega^q \le \sum_{j=1}^m \int_{p_j^*\{u_{jk} \not = v_{jk}\}} R_r \wedge \omega^q \lesssim \sum_{j=1}^m\|T_{jr}\|_{\{u_{jk} \not = v_{jk}\}} \le m \epsilon$$
and    similarly,
$$\int_{X^{m+1}} |\psi_k - \tilde{\psi}_k| R \wedge \omega^q \lesssim \epsilon$$
This combined with the fact that $\tilde{\psi}_k R_r \to \tilde{\psi}_k R$ as $r \to \infty$ implies that 
\begin{align}\label{limit-RrpaiktiendenRkjvofing}
\psi_k R_r \to \psi_k R
\end{align}
as $r \to \infty$.   Using this and the equality  $\psi_k+1 =0$ on $\cup_{j=1}^m p_j^*\{u_j \le -k\}$, we have 
\begin{align} \label{eq-bieudienpsikR}
-\psi_k R  &=- (\psi_k+1) R + R= - \lim_{r \to \infty}(\psi_k+1)R_r + R\\
\nonumber
&= - \lim_{r \to \infty}(\psi_k+1)R_k + R= R - (\psi_k +1)R_k.
\end{align}
Let $R_{k,\infty}$ be a tangent current to $-\psi_k R$ along $\Delta_{m+1}$. Note that $R_{k,\infty}$ is a positive current.  Since $\psi_k$ is bounded, we can apply Lemma \ref{le-hoitutoidensityboundedpsh} to $(\psi_k+1)R_k$. Consequently, 
\begin{align}\label{eq-pro-vanishingrestrictednonpluripolarRinftyTjk}
R_{k, \infty}= R_\infty  - \pi^* \bigg((1+ \rho_k ) \wedge_{j=1}^m T_{jk} \wedge T\bigg),
\end{align}
where $\rho_k:= \psi_k|_{\Delta_{m+1}}$.  Now observe
\begin{align*}
(1+ \rho_k) \wedge_{j=1}^m T_{jk} \wedge T &= (1+ \rho_k) \bold{1}_{\cap_{j=1}^m \{u_j > -k\}} \wedge_{j=1}^m T_{jk} \wedge T \\
&=  (1+ \rho_k) \bold{1}_{\cap_{j=1}^m \{u_j > -k\}} \langle \wedge_{j=1}^m T_j  \dot{\wedge}  T \rangle=(1+ \rho_k) \langle \wedge_{j=1}^m T_j  \dot{\wedge}  T \rangle 
\end{align*}
which converges to $\langle \wedge_{j=1}^m T_j  \dot{\wedge}  T \rangle$ as $k \to \infty$. This together with (\ref{eq-pro-vanishingrestrictednonpluripolarRinftyTjk}) yields the desired assertion and  finishes the proof. 
\endproof

For every closed positive current $S$ of bi-degree $(1,1)$ on $X$,  let $I_S$ be the polar set of potentials of $S$, \emph{i.e,} the set of points in $X$ where potentials of $S$ are equal to $-\infty$.  For $x \in X$, we denote $\nu(S, x)$ the Lelong number of $S$ at $x$. The following result allows one to estimate the h-dimension of density currents. 

\begin{proposition} \label{pro-uocluongmassonsmallsetdensity} %Let $T_1, \ldots, T_m$ be closed positive current of bi-degree $(1,1)$ and $T$ a closed positive current of bi-degree $(p,p)$ on a complex manifold.  
Let $A$ be a Borel subset of $X$. Assume that for every $J \subset \{1, \ldots, m\}$, the currents $(T_j)_{j \in J},T$ admit density currents and if $R_J$ is a density current associated to $(T_j)_{j \in J},T$ and $m_J:=|J|<m$, then   $R_J$ has  no mass on the set $$\pi_{m_J+1}^*(\cap_{j \not \in J} \{x \in A: \nu(T_j, x) >0\}).$$
 Then for every density current $S$ associated to $T_1, \ldots, T_m,T$, the h-dimension of the current $\bold{1}_A S$ is equal to $n-m-p$. In particular, if $p+m >n$, the current $\bold{1}_{X \backslash \cup_{j=1}^m I_{T_j}} S$ is zero, and if $p+m \le n$, the last current is  of minimal h-dimension.   
\end{proposition}

Note that $R_J$ is a current on $E_{m_J+1}$ and $\pi_{m_J+1}$ is the projection from  $E_{m_J+1}$ to $\Delta_{m_J+1}$. 

\proof  %The problem is local. Hence we can assume $X$ is K\"ahler and
Let $\omega$ is a Hermitian metric on $X$. Recall that we frequently identify $\Delta_{m+1}$ with $X$.  Let $q>n-m-p$ be an integer.  We need to show that $S \wedge \pi_{m+1}^* \omega^q$ has no mass on $A$. %Without loss of generality, we can assume that $A$ is relatively compact. 
 Let $\Phi$ be a compactly supported positive smooth form of bi-dimension $(q+ m+p, q+m+p)$ in the normal bundle $E_{m+1}$ of $\Delta_{m+1}$ in $X^{m+1}$. %Without loss of generality, we can assume that $A$ is relatively compact.  Fix a constant $\epsilon >0$. Let $U$ be a relatively compact open subset containing $A$ such that  $\bold{1}_{U \backslash A} S$ is of mass $\le \epsilom$. 

Let $S$ be a density current associated to $T_1, \ldots, T_m,T$ with a defining sequence $(\lambda_l)_l$. By extracting a subsequence of $(\lambda_l)_l$ if necessary, we can assume that  for every non-empty subset  $J \subset \{1,\ldots, m\}$,  there is a density current $R_J$  associated to $(T_j)_{j \in J}, T$ such that $(\lambda_l)_l$ is a defining sequence of $R_J$.  For every $J \in \{1, \ldots, m\}$,  let $p_J: X^{m+1} \to X^{m_J+1}$ be given by 
$$p_J(x_1, \ldots, x_{m+1})= \big((x_j)_{j \in J}, x_{m+1}\big).$$
The map $p_J$ induces a vector bundle map $\tilde{p}_J$ from $E_{m+1}$ to $E_{m_J+1}$.  Let  $\Omega_{J}$ be a compactly supported positive form of bi-dimension $(m_J+p, m_J+p)$ in $E_{m_J+1}$ such that $\Omega_J$ is strictly positive on an open neighborhood of  $\tilde{p}_J(\supp \Phi)$ in  $E_{m_J+1}$. 

Observe that $(\pi_{m+1})_*(S \wedge \omega^q \wedge \Phi)$ and $(\pi_{m_J+1})_*(R_J \wedge \Omega_{J})$ are measures on $X$.  Arguing as in the proof of \cite[Theorem 3.1]{VietTuanLucas} and using the fact that $q>n-m-p$,  we see that 
\begin{align}\label{ine-uocluonghdimensionSJ}
 (\pi_{m+1})_*(S \wedge \omega^q \wedge \Phi) \le  C  \sum_{J: \, |J|<m} \bigg(\prod_{j \not \in J} \nu(T_j, \cdot )\bigg) (\pi_{m_J+1})_*(R_J \wedge \Omega_{J})
%\langle  (\pi_{m+1})_*(S \wedge \omega^q \wedge \Phi), f \rangle  \le  C \sum_{J}\big \langle   (\pi_{m_J+1})_*(R_J \wedge \Omega_{J}),   f  \prod_{j \not \in J} \nu(T_j,\cdot) \big \rangle_{K_{r+1}},
\end{align}
for some constant $C$ independent of $T_1, \ldots, T_m,T$. The desired assertion immediately follows because the right-hand side of  (\ref{ine-uocluonghdimensionSJ}) has no mass on $A$ by hypothesis.  %For readers' convenience, we will briefly recall how to prove (\ref{ine-uocluonghdimensionSJ}). ???????????????????  Let $f$ be a continuous function compactly supported on $X$. $$  (\pi_{m+1})_*(S \wedge \omega^q \wedge \Phi) \le  C  \sum_{J: \, |J|<m} \big(\prod_{j \not \in J} \nu(T_j, \cdot )\big) (\pi_{|J|+1})_*(R_J \wedge \Omega_{|J|+1})$$ The measure on the right-hand side has no mass on $A$ by the hypothesis. Thus $ (\pi_{m+1})_*(S \wedge \omega^q \wedge \Phi)$ has no mass on $A$. 
This finishes the proof.
\endproof

\begin{proof}[End of the proof of Theorem \ref{the-main1}] Notice now that we work with a compact K\"ahler manifold $X$.  Assume first that  $\langle \wedge_{j=1}^m T_j \dot{\wedge} T \rangle = \curlywedge_{j=1}^m T_j \curlywedge T$. It follows that the class of  $\langle \wedge_{j=1}^m T_j \dot{\wedge} T \rangle$ is equal to that of $\wedge_{j=1}^m \{T_j\} \wedge \{T\}$ by Corollary \ref{cor-classDSproduct}. In other words, $T_1,\ldots, T_m$ are of $T$-relative full mass intersection.  We now prove the converse statement. 

Consider the case where $T_1,\ldots, T_m$ are of $T$-relative full mass intersection. We have 
$$\{\langle \wedge_{j=1}^m T_j \wedge T\}= \wedge_{j=1}^m \{T_j\} \wedge \{T\}.$$
  By \cite[Lemma 4.10]{Viet-generalized-nonpluri}, for every $J\subset \{1,\ldots,m\}$, the currents $(T_j)_{j \in J}$ are of $T$-relative full mass intersection and  $T$ has no mass on $I_{T_j}$.  

We prove the desired assertion by induction on $m$.  When $m=0$, it is clear. We assume it holds for every $m'<m$. We prove it for $m$. By induction hypothesis,  for every $J \subsetneq \{1, \ldots,m\}$, the Dinh-Sibony product  of $(T_j)_{j \in J}, T$ is well-defined and equal to $\langle \wedge_{j \in J} T_j \dot{\wedge} T \rangle $. The last current  has no mass on $I_{T_l}$ because $T$ does so,  for $1 \le  l \le m$.  This allows one to apply Proposition \ref{pro-uocluongmassonsmallsetdensity} to $A:= X$, $T_1, \ldots, T_m,T$ to see that the density h-dimension of $T_1, \ldots, T_m,T$ is minimal. We deduce that   
$$\pi_{m+1}^* \big(\wedge_{j=1}^m \{T_j\} \wedge \{T\}\big)$$ (which is equal to the pull-back by $\pi_{m+1}$ of the class of $\langle \wedge_{j=1}^m T_j \wedge T\rangle$)
 is the total density class of $T_1, \ldots, T_m,T$. This combined with Theorem \ref{the-phanbucuarestricteddenstyvanonpluri}  gives the desired assertion.  The proof is finished.
\end{proof}

\bibliography{biblio_family_MA,biblio_Viet_papers}
\bibliographystyle{siam}

\bigskip

\noindent
\Addresses
\end{document}